\magnification=\magstep1
\centerline {\bf DEVELOPING COMPUTER PROGRAMS FOR KNOT CLASSIFICATION}
\bigskip
\centerline {Charilaos Aneziris}
\centerline {\it Institut f\"ur Hochenergienphysik Zeuthen}
\centerline {\it DESY Deutsches Elektronen-Synchrotron}
\centerline{\it Platanenallee 6}
\centerline {\it 15738 Zeuthen}
\bigskip
\centerline {\bf Abstract}
\bigskip
In this paper we summarise the work discussed in Ref. [1] and [2], in which we
introduced a method helpful in
solving the problem of knot classification. We also  
present results obtained since then.
\bigskip
\centerline {1. \qquad INTRODUCTION}
\bigskip
Knot Theory has attracted significant attention during recent years, both 
among mathematicians, and among areas of applied science such as Physics, 
Chemistry and Biology. In fact, a number of problems that were previously
considered unrelated to each other, have been connected through applications of
Kont Theory. While enormous progress has been achieved in the study 
of knots and of their applications, the problem of a complete classification
remains still open, in spite of recent successes (Ref. [3]). In this paper we
describe and discuss an algorithmic approach that could be useful in solving
the problem. With the help of an algorithm which is presented in this paper, a 
computer program was developed, resulting to the classification of all knots
with crossing number up to 11.
\medskip
In Section 2 we present the main points of the algorithm, while in Section 3 we
introduce a suitable notation and show how through this notation it is possible 
to classify knot projections. In Section 4 we show how Reidemeister moves can 
be used to identify projections of equivalent knots, so that ambient isotopic 
knots
may not appear more than once at the output. In Section 5 we generate a series
of ``color tests" in order to demonstrate knot inequivalence; such a procedure
is necessary since equivalent knots may fail to be
identified through the procedure of Section 4. Finally in Section 6 we show the
results obtained through this computer program. Ideally, any two knot 
projections should either belong to knots shown equivalent  
through Reidemeister moves, or to knots shown inequivalent due to different
responses in one or more ``color tests". This would be the case if the computer
program could run for ever; in practice the results depend on the two input
parameters, one indicating the maximum crossing number considered, the other
indicating the ultimate ``color test" to be used. Currently this has been
achieved for all knots whose crossing number does not exceed 11.
\bigskip
\centerline {2. \qquad THE ALGORITHMIC PROCESS}  
\bigskip
In this Section we present the main steps of the algorithm. First, a suitable
method to denote knot projections is introduced. Second, once the set of
possible such notations has been obtained, one needs a method to distinguish
notations that correspond to actual knot projections, from notations that do 
not. Third,         
notations that correspond to identical knot projections must be identified. 
Once these
steps are completed, knot projections are fully classified. This however is not
identical to classifying knots, since distinct knot projections may correspond
to ambient isotopic knots, and such knots are considered equivalent.
\medskip
The next step therefore is to identify such projections. It is well known that
projections of ambient isotopic knots are related through Reidemeister moves
(Ref. [4]). It is thus necessary to know how a notation is affected by a
Reidemeister move. Once this is known, one may use such moves to identify
ambient isotopic knots. In order for the program to be finite, one may
establish an upper limit to the number of Reidemeister moves to be applied; it
turns out however to be simpler to set an upper limit 
to the crossing number of the knot projections involved, instead of the
number of Reidemeister moves
This upper limit is one of the two input parameters used in the program.
\medskip
Since however no upper limit to either this number, or the number of
necessary moves is known, there is no certainty that projections not found 
connected through Reidemeister moves, will actually belong to inequivalent 
knots.
Therefore one starts by identifying as many equivalent knots as possible.  
Then one 
proceeds by selecting one knot from each equivalence class, conventionally the
knot appearing first, and by calculating knot characteristics, in order to
establish inequivalences between selected knots. As such characteristics we
shall use the so called ``color tests", which are a generalisation of the
``tricolorisation" through which the trefoil's non-triviality may easily be
shown. Each color test consists of an $n \times n$ matrix whose elements take
values in $\{ 1,2,...,n \}$; the strands of each knot projections are mapped to
elements of $\{ 1,2,...,n \}$ (the $n$ ``colors"). Acceptable mappings are the 
ones where the three strands meeting at each crossing, are mapped to numbers
satisfying relations determined through the $n \times n$ matrix.
Once certain constraints among the
matrix elements are satisfied, the number of acceptable mappings is invariant
under Reidemeister moves. Therefore if two projections yield different results
for one or more such color tests, they definitely belong to inequivalent knots.
\medskip
Not all knots on which such ``color tests" are applied, are going to yield
distinct results and thus shown inequivalent. This is due to two reasons.
First, some of these knots are actually equivalent, but due to the limitations
in the Reidemeister moves considered, the program failed to identify them. 
Second,
even if two knots are actually inequivalent, they may not yield distinct 
results due 
to the finite number of color tests applied. The second input parameter 
indicates in fact the color tests that are applied.
\medskip
Having presented the main steps of the program, we now proceed with a detailed 
discussion. 
\bigskip
\centerline {3. \qquad THE CLASSIFICATION OF KNOT PROJECTIONS}
\bigskip
Knot projections are denoted as sets of $n$ pairs of natural numbers \break 
$\{ (a_1,a_2),(a_3,a_4),...,(a_{2n-1},a_{2n}) \}$, where $n$ is the crossing
number, such that \break 
$a_i \in \{1,2,...,2n\}$ and $i \neq j \Leftrightarrow
a_i \neq a_j$. This set is obtained as follows. First one chooses a starting
p oint and an orientation. Then, as one travels around the projection, one
assigns successive natural numbers to the crossing points, starting from $1$
and ending to $2n$. Each crossing is eventually assigned two numbers, 
$a_{over}$ for the overcrossing and $a_{under}$ for the undercrossing. The set
of the pairs $(a_{over},a_{under})$ denotes the projection.
\medskip
Not all possible notations yield actual knot projections, the simplest
counterexample being $\{(1,3),(2,4)\}$. One necessary condition is that odd
numbers are always paired to even numbers. This condition is not sufficient, as
the counterexample \break
$\{(1,4),(3,6),(5,8),(7,10),(9,2)\}$ demonstrates. The
necessary and sufficient condition is that any two loops obtained from an
actual projection, must either share one or more line segments, or intersect at
an even number of points, vertices not being counted. This condition is due to
the {\it Jordan Curve Theorem} (Ref. [5]) which states that any loop on $R^2$
or $S^2$ which does not intersect itself, divides $R^2$ or $S^2$ into two
disjoint pieces. In these two counterexamples, the loops 1-2-3 and 3-4-1, in
the first case, and 1-2-3-4 and 5-6-7-8 in the second case, violate this rule
by not sharing any common segment and intersecting at exactly one point. 
The maximum number of loops obtained from an $n$ crossing knot
projection is $3^n$, since each crossing may be a vertex of the loop or may 
not,
and if it is, there are two possible direction changes. Therefore, checking
whether a notation yields an actual knot projection, is a finite process.  
\medskip
Once ``drawable" notations have been separated from ``undrawable" ones, one
needs to identify notations leading to the same knot projection. For an $n$
crossing projection, there are at most $4n$ such projections, corresponding to 
$2n$ possible starting segments and to $2$ possible orientations. By altering 
the starting place and/or the orientation, each pair $(a_i,a_j)$ becomes 
$(k+\epsilon a_i,k+\epsilon a_j)$, where $k$ indicates the change of the 
starting
point and $\epsilon =\pm 1$ indicates a possible change of orientation; 
$\epsilon = 1$ indicates that the orientation remains the same, while
$\epsilon = -1$ indicates that the orientation has been reversed. One may thus 
identify all such notations and keep just one, conventionally the one appearing
first. This too is a finite process, and since it is less time consuming than
checking the notation's ``drawability", the program becomes more efficient if
this step precedes the previous one.
\medskip
At this point the procedure of classifying two-dimensional knot projections has
been completed.
\bigskip
\centerline {4. \qquad IDENTIFYING AMBIENT ISOTOPIC KNOTS}
\bigskip
As mentioned earlier, two projections correspond to ambient isotopic knots if
and only if they can be connected through Reidemeister moves. There are three
kinds of Reidemeister moves, and their pictorial forms can be found in a
number of relevant books (see for example Ref. [6]). Here we present their
``numerical" form, by showing how each Reidemeister move affects a 
notation.
\medskip
A first Reidemeister move, which increases the crossing number by 1, adds a
pair $(i,i+1)$ or a pair $(i+1,i)$ to the notation, while replacing any other
number $j$ which is larger or equal to $2$, with $j+2$. A second Reidemeister
move, which increases the crossing number by 2, adds two pairs $(i,j)$ and
$(i+1,j+1)$, or $(i+1,j)$ and $(i,j+1)$, to the notation. Numbers larger or 
equal to $i$ and smaller than
$j$, increase by $2$; numbers larger or
equal to $j$ increase by $4$. A first or second Reidemeister move which
decrease the crossing number, will have the converse effect. Finally a third 
Reidemeister move, which keeps the crossing number constant, replaces pairs
$(i,j)$, $(i',k)$ and $(j',k')$ with the pairs $(i,k')$, $(i',j')$ and $(j,k)$,
where $|i'-i|=|j'-j|=|k'-k|=1$, while all other $n-3$ pairs remain the same.
\medskip
The process of identifying equivalent knot goes as follows. First, one obtains
through the procedure of Section 3, all distinct knot projections whose
crossing number does not exceed some maximum value $N$. Then, on each 
projection
one applies Reidemeister moves that do not increase the crossing number. 
Projections that cannot be connected to ones appeared before, are stored in
the computer memory and are assigned two numbers, a ``temporary" and
a ``permanent" one. Initially these numbers are equal. The permanent numbers
assigned to such projections, are successive natural numbers. Projections
connected to ones appeared before, are not stored in the memory, but help
obtain equivalences among projections already stored. If for example
some projection $P$ is found equivalent to projections $P_1$, $P_2$, ..., $P_k$
which have been assigned the permanent numbers $p_1$, $p_2$, ..., $p_k$ and the
temporary numbers $t_1$, $t_2$, ..., $t_k$, the permanent numbers do not 
change, while the temporary numbers are replaced by $\min (t_1,t_2,...,t_k)$.
\medskip
When all projections have been checked, only the ones stored in the memory
with equal temporary and permanent numbers are going to appear at the output,
since only these have not been found equivalent to preceding projections. As
stated earlier, such projections may or may not be equivalent, and one thus
procedes by developing ``color tests" in order to distinguish inequivalent
knots. 
\bigskip
\centerline {5. \qquad ESTABLISHING ``COLOR TESTS"}
\bigskip
A simple method to show the existence of non-trivial knots is through
``tricolorisation". One maps the strands $s_i$ of a knot projection to a
number $n_i \in \{ 1,2,3 \}$, so that at each crossing, the strands involved,
$s_i$, $s_{i+1}$ and $s_j$ satisfy the relation $n_i+n_{i+1}+n_j=0{\rm mod}3$.
If the projection is altered by a Reidemeister move, to each mapping of the
old projection corresponds exactly one mapping of the new. Therefore the
number of mappings is a knot invariant; if a projection $P_1$ admits $m_1$
mappings, while a projection $P_2$ admits $m_2$, and $m_1 \neq m_2$, then
$P_1$ and $P_2$ definitely belong to inequivalent knots. For the trefoil
three such mappings are possible, each mapping the only strand to one of the
elements of $\{ 1,2,3 \}$. In contrast, for the trefoil nine such mappings are
possible; three map all strands to the same number, while the other six map
the strands to three different numbers. Therefore the non-triviality of the
trefoil is established (Ref. [7]). 
\medskip
Starting from this ``three color test", one may generalise to obtain more such
color tests in order to distinguish inequivalent knots whose responses to
tricolorisation are identical. Each such color test is defined through an
$n \times n$ matrix $M_{ij}$, so that if at some crossing the strands involved,
$s_i$, $s_{i+1}$ and $s_j$ are mapped to $n_i$, $n_{i+1}$ and $n_j$, then
either $n_{i+1}=M_{n_in_j}$, or $n_i=M_{n_{i+1}n_j}$, depending on whether the
crossing is positive or negative. Only mappings $s_k \rightarrow n_k$, where
this property is satisfied at every crossing, are considerd acceptable and are
counted for the corresponding knot invariant. For the ``three color test"
mentioned before, one may notice that $n=3$, $M_{ii}=i$, while for $i \neq j$,
$M_{ij}=k$, where $k \neq i$ and $k \neq j$.
\medskip 
Not all possible matrices however are suitable. A matrix may only be used
to define a ``color test" if for any two knot projections $P$ and $P'$
differing by Reidemeister moves, to each acceptable mapping for $P$ 
corresponds exactly one mapping for $P'$. To ensure this property, one
considers the constraints that each Reidemeister move imposes. One may easily
observe that these constraints are the following.
\medskip 
$1^{st}$ move: $M_{ii}=i \quad \forall \quad i \in \{ 1,2,...,n \}$
\smallskip
$2^{nd}$ move: $M_{ij}=M_{i'j} \Leftrightarrow i=i'$
\smallskip
$3^{rd}$ move: $M_{ij}=k \quad \wedge \quad M_{li}=m \quad \wedge \quad 
M_{lj}=n
\quad \Rightarrow M_{nk}=M_{mj}$
\medskip In addition, an $n$-color test is not considered if there is a subset
$S$ of $\{ 1,2,...,n \}$ other than the empty set and $\{ 1,2,...,n \}$ itself,
such that $i \in S \Rightarrow M_{ij} \in S \quad \forall \quad j \quad \in
\quad \{ 1,2,...,n \}$, since such a test may be reduced to simpler ones.
Finally, two tests are considered identical if one may be obtained from the 
other by permutation, or if they are defined through matrices $M$ and $M'$
such that $M_{ij}=k \Rightarrow M'_{kj}=i$, since in such a case they are
related through mirror symmetry.
\medskip
Subsequently, a computer program was developed that recorded the matrices that 
yield
distinct valid color tests. The running time grew
exponentially with the number of colors; to obtain all color tests for up to 11
colors the time needed was a few days, while for 12 colors it would exceed one
month. The number of color tests per number of colors came out to be as
follows.
\bigskip \hskip 2 cm 
\vbox{\catcode`\*=\active \def*{\hphantom{0}}
\offinterlineskip
\halign{\strut#&\vrule#\quad&\hfil#\hfil&\quad\vrule#\quad&\hfil#\hfil&\quad
\vrule#\cr
\noalign{\hrule}
&&\omit {\bf Number of Colors} && {\bf Number of Tests} &\cr \noalign {\hrule}
&&*1&& *1&\cr 
&&*2&& *0&\cr 
&&*3&& *1&\cr 
&&*4&& *1&\cr 
&&*5&& *2&\cr 
&&*6&& *2&\cr 
&&*7&& *3&\cr 
&&*8&& *2&\cr 
&&*9&& *6&\cr 
&&10&& *1&\cr 
&&11&& *5&\cr \noalign{\hrule}}}   
\bigskip
As shown in Section 6, these tests are not sufficient for distinguishing knots
of high crossing numbers,
and the method of establishing tests by explicitly checking every possible
matrix is not efficient enough. Instead, one obtains an infinite number of 
tests by
generalising from the tests already established. One such class of tests is
defined through matrices \break $M_(i,j)=(k+1)j-ki \quad {\rm mod} \quad n$, 
where the
Greatest Common Divisor \break 
GCD($k,n$)=GCD($k+1,n$)=1. The existence of non-trivial
mappings depends on the determinant of the linear homogeneous system that is 
defined
through the equations satisfied at each crossing. This determinant is the
{\it Alexander-Conway} polynomial (Ref.[8]). One may thus calculate and compare 
the Alexander-Conway polynomials of various knots, and apply additional color
tests only for knots whose Alexander-Conway polynomials are identical.
\medskip
A second class of tests associates the ``colors" to group elements $g_i$, and
is defined through the matrix $M(g_i,g_j)=g_jg_ig_j^{-1}$ (Ref. [9]). In
particular, one may use as groups the permutation groups $S_n$; each 
conjugacy class, defined through a partition $\lambda _1 \geq \lambda _2 \geq
... \geq \lambda _k$ of $n$, ($\lambda _1 + \lambda _2 + ... + \lambda _k = 
n$), defines a valid color test.
\bigskip
\centerline {6. \qquad COMPUTER RESULTS}
\bigskip
The maximum value of the crossing number of the projections studied, was set
equal to $N=14$. In order for the program to run, the CPU time needed was 8
days, and the memory required was about 10 MBytes. The number of knots that
were not connected through Reidemeister moves, came out as follows.  
\bigskip \hskip 2 cm 
\vbox{\catcode`\*=\active \def*{\hphantom{0}}
\offinterlineskip
\halign{\strut#&\vrule#\quad&\hfil#\hfil&\quad\vrule#\quad&\hfil#\hfil&\quad
\vrule#\cr
\noalign{\hrule}
&&\omit {\bf Number of Crossings} && {\bf Number of Knots} &\cr \noalign 
{\hrule}
&&*0&& ****1&\cr 
&&*1&& ****0&\cr 
&&*2&& ****0&\cr 
&&*3&& ****1&\cr 
&&*4&& ****1&\cr 
&&*5&& ****2&\cr 
&&*6&& ****3&\cr 
&&*7&& ****7&\cr 
&&*8&& ***21&\cr 
&&*9&& ***49&\cr 
&&10&& **165&\cr 
&&11&& **552&\cr 
&&12&& *2191&\cr 
&&13&& 29781&\cr \noalign{\hrule}}}   
\bigskip
Due to memory constraints, 14 crossing knots were not recorded. As pointed out
earlier, these numbers are mere upper limits, since it is certain that many of
them although equivalent, may only be connected through Reidemeister moves
involving more than 14 crossings. One thus proceeds by applying the color tests
in order to obtain topologically inequivalent knots.
\medskip
When the color tests listed in the table of Section 5 were applied, which are
all the 
color tests involving at most 11 colors, all knots with crossing number up to 7
were shown inequivalent. This was not the case however with knots whose
crossing number is 8, and therefore this method is good enough for only the
first 15 knots.
\medskip
When the Alexander-Conway polynomials were calculated, the results were
slightly better; all 36 knots whose crossing number does not exceed 8, possess
distinct Alexander-Conway polynomials. When knots with crossing number 9 are
also considered, one faces the first cases of inequivalent knots with identical
Alexander-Conway polynomials.
\medskip
We later applied color tests derived from permutation groups, as discussed at
the end of Section 5. Permutation groups up to $S_5$ were sufficient to
demonstrate the inequivalence of all knots that possess identical 
Alexander polynomials and whose crossing number does not exceed 10. For 
crossing number 11, one has to go up to $S_7$, until all 802 knots with
crossing number not exceeding 11 were shown inequivalent. For a complete list 
of all these 
knots and the characteristics through which these knots were distinguished, the
reader is referred to Ref. 10.
\medskip
For crossing numbers 12 and 13 it is almost certain that equivalent knots do
exist, which would require the study of projections with crossing number 
higher than 14. The basis of this assumption is the fact that 11 crossing knots
may only be distinguished once 14 crossing projections are studied; if the
maximum value is set equal to 13 crossings, then 3 pairs of 11 crossing ambient
isotopic knots cannot be identified.
\medskip
At this point we have derived a full list of all knots whose crossing number 
does not exceed 11. In principle the method discussed could lead to extending
this list to an arbitrary high crossing number; the CPU time and computer
memory however rise very rapidly with the crossing number. 
\bigskip
\centerline {A C K N O W L E G M E N T S}
\bigskip
I would like to thank the Brookhaven National Laboratory, the Weizmann 
Institute of Science and the University of G\"ottingen for their kind 
hospitality during parts of the work. I would also like to thank Wei Chen,
K. Anagnostopoulos, J. Przytycki, S. Lambropoulou, B. Westbury, C. Shubert,
S. Garoufalidis for
their help and/or interesting discussions concerning the problem discussed in
this paper.   
\bigskip
\centerline {R E F E R E N C E S}
\bigskip
1) C. Aneziris, DESY-Ifh preprint DESY-94-230, November 1994.
\par \smallskip
2) C. Aneziris, q-alg/9505003, May 1995.
\par \smallskip
3) M.B. Thistlethwaite, (1985) L.M.S. Lecture Notes no 93 pp. 1-76, Cambridge
University Press.
\par \smallskip
4) K. Reidemeister, {\it Abh. Math. Sem. Univ. Hamburg} {\bf 5} (1927) 7-23.
\par \smallskip
5) C. Jordan, (1893) {\it Course d' Analyse}, Paris.
\par \smallskip
6) D. Rolfsen, (1976) {\it Knots and Links}, Berkeley, CA: Publish or Perish,
Inc., L.H. Kauffman, (1991) {\it Knots and Physics}, World Scientific
Publishing Co. Pte. Ltd. and references therein.
\par \smallskip 
7) R.H. Fox, {\it Canadian J. Math.}, XXII(2) 1970, 193-201. 
\par \smallskip 
8) J.W. Alexander, {\it Trans. Amer. Math. Soc.} {\bf 30} (1926) 275-306.
\par \smallskip 
9) R. Fenn and C. Rourke, {\it Journal of Knots and its Ramifications},
1(4) 1992, 343-406.
\par \smallskip 
10) C. Aneziris, http://sgi.ifh.de/\~{}aneziris/contents/
\end